\documentclass[epjST]{svjour}
\usepackage{amsmath}
\usepackage{amssymb}
\usepackage{graphics}
\usepackage{epsfig}
\usepackage{epstopdf}
\usepackage{algorithm}
\usepackage{algpseudocode}
\usepackage{enumitem}
\usepackage{tabularx}
\usepackage{caption}
\usepackage{authblk}
\usepackage{wrapfig}
\usepackage{cite}
\usepackage{pdfpages}
\usepackage[figuresright]{rotating}
\usepackage[toc,page]{appendix}

\def\bs{\boldsymbol}
\newcommand{\rom}[1]{\uppercase\expandafter{\romannumeral #1\relax}}
\begin{document}
\title{Identifying manifolds underlying group motion in Vicsek agents}
\author{Kelum Gajamannage\inst{1}\fnmsep\thanks{\email{dineshhk@clarkson.edu}} 
\and Sachit Butail\inst{2}\fnmsep\thanks{\email{sbutail@iiitd.ac.in}}
\and Maurizio Porfiri\inst{3}\fnmsep\thanks{\email{mporfiri@nyu.edu}}
\and Erik M. Bollt\inst{1}\fnmsep\thanks{\email{ebollt@clarkson.edu}}}
\institute{Department of Mathematics, Clarkson University, Potsdam, New York, USA. 
\and Indraprastha Institute of Information Technology Delhi (IIITD), New Delhi, India. 
\and Department of Mechanical and Aerospace Engineering, New York University, Brooklyn, New York, USA.}
\abstract{
Collective motion of animal groups often undergoes changes due to perturbations. In a topological sense, we describe these changes as switching between low-dimensional embedding manifolds underlying a group of evolving agents. To characterize such manifolds, first we introduce a simple mapping of agents between time-steps. Then, we construct a novel metric which is susceptible to variations in the collective motion, thus revealing distinct underlying manifolds. The method is validated through three sample scenarios simulated using a Vicsek model, namely switching of speed, coordination, and structure of a group. Combined with a dimensionality reduction technique that is used to infer the dimensionality of the embedding manifold, this approach provides an effective model-free framework for the analysis of collective behavior across animal species.
} %end of abstract
\maketitle

%%%%%            Section 01:	Introduction          %%%%%
\section{Introduction}\label{sec:introduction}
In animal groups, the response to a perturbation---internal or external---is often manifested in the form of changes in group speed, coordination, or structure \cite{carere2009aerial, couzin2002collective, helbing2000simulating, neill1974experiments, vicsek2012collective}. Such changes are witnessed in fish schools and bird flocks under attack \cite{vabo1997individual, parrish1991predators, lima1990behavioral}, foraging animal groups \cite{deneubourg1989collective, couzin2009collective}, and human crowds exposed to alarm situations leading to panic \cite{helbing2002simulation, shiwakoti2011animal}. Based on our recent effort demonstrating that collective motion is associated with a low-dimensional embedding \cite{abaid2012topological, delellis2014collective, delellis2013topological, butail2013analysis, gajamannage2015dimensionality}, we expect that such behavioral changes should be manifested in variation of the topology of an underlying manifold.

Coordinated group behaviors can be represented as a manifold, $\mathcal{M}$, in an abstract phase space of a group such that configurations $\bs{\mathcal{A}}$ evolve in time, $t$ according to the mapping
\begin{equation}\label{eqn:mapping_manifold}
\Phi^{(t)}(\bs{\mathcal{A}}):{\cal M}\rightarrow {\cal M}.
\end{equation}
The presence of a phase change, could therefore be defined as a union of two distinct underlying manifolds $\mathcal{M}^{(j)}, j=1 \ \text{and} \ 2$. 

\begin{figure}[htp]
        	\centering
	\vspace{-10pt}
        	\includegraphics[width=1\textwidth]{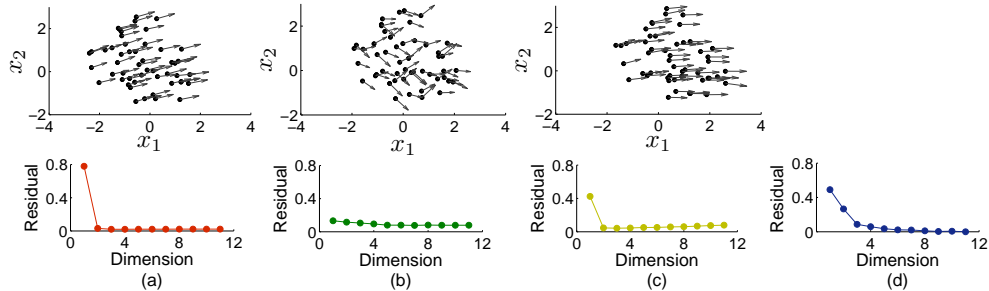}
	\vspace{-10pt}
        	\caption{A simulation of self-propelled particles in two dimensions. Changes in coordination can be seen as switching of embedding manifolds of different dimensions. The top row shows snapshots of self-propelled particles as they move with low, high, and then low noise again. The bottom row displays the normalized Isomap residual variance  for each instance (a, b, c) and the whole motion (d). The location of elbow in each plot indicates that  low noise instances are low-dimensional (3), high noise instance results into a high-dimensional embedding (6). The embedding dimensionality for the whole motion is four.
\vspace{-10pt}
}\label{fig:vicsek_reversion}
\end{figure}

As an illustrative example of change between high and low coordination in a group, we simulate 20 self-propelled particles using the Vicsek model \cite{vicsek1995novel} by alternatively imposing low and high noise to their individual dynamics in three distinct phases. Figure \ref{fig:vicsek_reversion} shows the resulting variation in the coordination of agents. Following our previous work \cite {abaid2012topological, delellis2014collective, delellis2013topological, butail2013analysis}, we run an established dimensionality reduction routine called Isomap \cite{tenenbaum2000global} on configurations from each phase as well as the whole dataset comprising all three phases to infer the dimensionality of the underlying embedding. The residual of the reconstruction error shows that the first and last instances (low noise) embed on manifolds with lower dimensionality than the middle instance (high noise) \cite{abaid2012topological}. However, when the same routine is applied to the full dataset, we find an embedding dimensionality that is in neither an indicator of high nor low group coordination. This example suggests that group behavior is manifested in distinct embedding manifolds, and a naive implementation of dimensionality reduction could result in a loss of information about critical changes in collective behavior \cite{abaid2012topological}.

In this paper, we propose a method to study higher-dimensional alternating manifolds in response to an internal or external perturbation to a group of agents. Our approach shares similarities with the coarse-grained analysis of stochasticity-induced switching in one-dimensional model of collective motion \cite{kolpas2007coarse}. Therein, a single coarse observable, the average nearest-neighbor distance, is used to distinguish between the stationary and mobile states in the collective motion.

Here, we characterize group motion by using a single coarse observable $\Delta$, which computes the distance between sub-manifolds corresponding to distinct group actions. Specifically, we use agent positions in each configuration to create a single weighted metric of group speed, orientation, and structure. The positions of agents between consecutive time-steps are mapped using a nearest-neighbor search in position and velocity, which suffices to characterize the group motion. We succinctly present our approach in Section \ref{sec:method}, with further technical details on the velocity computation reported in the Supplementary materials. In Section \ref{sec:examples}, we assess the performance of our approach on three synthetic examples of collective behavior using the Vicsek model \cite{vicsek1995novel} simulating changes in group speed, coordination, and structure. We conclude in Section \ref{sec:conclusion} with a discussion of the performance of the method and ongoing work.

%%%%         Section 02: 	Method of characterizing manifolds       %%%%
\section{A method for characterizing alternating manifolds}\label{sec:method}
We assume that agents move in a two-dimensional space, so that the raw data for our method  consists of two-dimensional position vectors of all the agents in time. For a group of $N$ agents, the position vector of the $i$-th agent at the $t$-th time step is denoted by $\bs{a}^{(t)}_i\in \mathbb{R}^2$. The configuration vector, defined as
\begin{equation}\label{eqn:configuration}
\bs{\mathcal{A}}^{(t)}=[\bs{a}^{(t)}_1; \bs{a}^{(t)}_2; \dots; \bs{a}^{(t)}_N]\in \mathbb{R}^{2N},
\end{equation}
is a point on the manifold that represents the position of the entire group at time $t$. The order of the agents in the configuration vector is not required to be preserved for this method, since the mapping described herein determines the order of each agent in consecutive configurations. The data acquired from the group at distinct times $t=1, 2, \dots, T$ is consolidated in the dataset $\mathcal{Q}=\{\bs{\mathcal{A}}^{(t)}\}_{t=1}^T$. If the complete dataset $\mathcal{Q}$ of configurations is partitioned into $m$ distinct sub-manifolds, such that ${\cal M}=\cup_{j=1, ..., m}   {\cal M}^{(j)}$, the dimensionality of each sub-manifold can then be found by implementing a dimensionality reduction algorithm, such as Isomap \cite{tenenbaum2000global}.

We describe the underlying manifold using the speed, coordination, and structure of the agents in the group. In order to compute group speed, we calculate the velocity of each from a time step to the next one. For that, we first construct a bijective evolution function, $\Phi^{(t)}$ in (\ref{eqn:mapping_manifold}), for agents between time steps. The identity of agents in the input data is not required to be preserved for this method, as the mapping presented below automatically finds any agent's position change between consecutive time steps. Given that this mapping has a rather modest computational cost of $O(\log N)$, we opt for this approach in favor of more accurate but computationally costly methods \cite{khan2005mcmc,  wu2009tracking-a}. 

We assume that the position of an agent at the next step are in a neighborhood of its position at the previous time step.
\begin{figure}
        	\centering
        	\includegraphics[width=.4\textwidth]{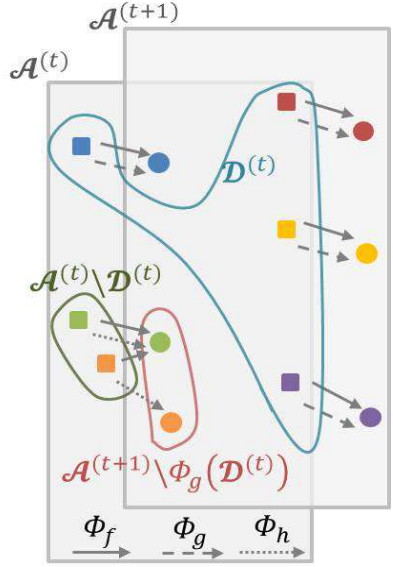}
        	\caption{Mapping of positions of agents between two consecutive time-steps. Same colored square and circle pair refers to the same agent. Squares represent agents' positions at the $t$-th time step and circles represent agents' positions at the $(t+1)$-th time-step.}
	\label{fig:maping}
\end{figure}
We search the nearest neighbor \cite{friedman1977algorithm}, say $\bs{a}_j^{(t+1)}\in \bs{\mathcal{A}}^{(t+1)}$, for an agent, say $\bs{a}^{(t)}_i\in \bs{\mathcal{A}}^{(t)}$, and map them through $\Phi^{(t)}_f$ as
\begin{equation}\label{eqn:mapping}
\Phi^{(t)}_f: \bs{\mathcal{A}}^{(t)} \rightarrow \bs{\mathcal{A}}^{(t+1)} \ \text{such that} \ \Phi^{(t)}_f \left(\bs{a}_i^{(t)}\right):=\bs{a}_j^{(t+1)}.
\end{equation}
Similarly, we apply $\Phi^{(t)}_f$ to all the agents at the $t$-th time step to map the corresponding agents at the $(t+1)$-th time step. Since some agents at the $t$-th time step may share the same nearest neighbor at $(t+1)$-th time step, $\Phi^{(t)}_f$ maps more than one agents in $\bs{\mathcal{A}}^{(t)}$ into the same agent in $\bs{\mathcal{A}}^{(t+1)}$ and violate the uniqueness of the mapping. In order to ensure the uniqueness of the mapping, $\Phi^{(t)}_f$ should be a bijection which is defined to be a one-to-one correspondence of agents between consecutive time steps. We define and extract the sub-domain $\bs{\mathcal{D}}^{(t)} \subseteq \bs{\mathcal{A}}^{(t)}$ as all agents bijectively mapped with agents in $\bs{\mathcal{A}}^{(t+1)}$ through $\Phi^{(t)}_f$. The rest of the agents are mapped using a different map that assures the bijection, as described in what follows. We denote the bijective mapping $\Phi^{(t)}_f$ in $\bs{\mathcal{D}}^{(t)}$ by $\Phi^{(t)}_g$
\begin{equation}\label{eqn:mapping_f}
\Phi^{(t)}_g: \bs{\mathcal{D}}^{(t)} \rightarrow \bs{\mathcal{A}}^{(t+1)}.
\end{equation}
Mappings $\Phi^{(t)}$ and $\Phi^{(t)}_g$ and the domain $\bs{\mathcal{D}}^{(t)}$ for six agents are illustrated in Figure \ref{fig:maping}. For convenience and without lack of generality, we assume that the first $m$ agents at the $t$-th time step are in $\bs{\mathcal{D}}^{(t)}$. The velocity of the $i$-th agent in $\bs{\mathcal{D}}^{(t)}$, say $\bs{v}^{(t)}_i \in \mathbb{R}^2$ for $i=1,2, \dots m$, is defined as the displacement from the current positions to the next position
\begin{equation}\label{eqn:velocity_f}
\bs{v}^{(t)}_i=\Phi^{(t)}_g(\bs{a}^{(t)}_i)-\bs{a}^{(t)}_i.
\end{equation} 
We compute the mean velocity of the agents in $\bs{\mathcal{D}}^{(t)}$ as $\bs{\mu}^{(t)}_1= \text{mean}\{\bs{v}^{(t)}_i \ ; \ i=1,2, \dots, m\}$ and use this quantity to map the remaining elements in $\bs{\mathcal{A}}^{(t)}$.

Mapping of the remaining elements in the configuration $\bs{\mathcal{A}}^{(t)}$, denoted as $\bs{\mathcal{A}}^{(t)} \backslash \bs{\mathcal{D}}^{(t)}$, with the remaining elements in $\bs{\mathcal{A}}^{(t+1)}$, denoted as $\bs{\mathcal{A}}^{(t+1)} \backslash \Phi^{(t)}_g (\bs{\mathcal{D}}^{(t)})$, is presented below. We choose an agent $\bs{a}_i^{(t)}\in \bs{\mathcal{A}}^{(t)}\backslash \bs{\mathcal{D}}^{(t)}$ for $i=m+1, \dots, N$ and compute displacement components from that agent to all the agents $\bs{a}_j^{(t+1)} \in \bs{\mathcal{A}}^{(t+1)}\backslash \Phi^{(t)}_g(\bs{\mathcal{D}}^{(t)});$ $j=m+1,\dots, N$. We use the mean velocity, $\bs{\mu}^{(t)}_1$, to approximate the positions of unmapped agents in the $t$-th time step. If the displacement computed between $\bs{a}_i^{(t)}$ and $\bs{a}_j^{(t+1)}$ is the closest to the mean velocity of the agents mapped bijectively, $\bs{\mu}^{(t)}_1$, then $\bs{a}_i^{(t)}$ is mapped to $\bs{a}_j^{(t+1)}$. Thus, we map all the agents bijectively between time steps $t$ and $t+1$ such that they all follow the same group action. We apply this mapping for all the agents in $\bs{\mathcal{A}}^{(t)}\backslash \bs{\mathcal{D}}^{(t)}$.  The mapping between agents $\bs{\mathcal{A}}^{(t)} \backslash \bs{\mathcal{D}}^{(t)}$ and $\bs{\mathcal{A}}^{(t+1)} \backslash \Phi^{(t)}_g (\bs{\mathcal{D}}^{(t)})$ is defined as 
\begin{equation}\label{eqn:mapping_g}
\Phi^{(t)}_h:\bs{\mathcal{A}}^{(t)}\backslash \bs{\mathcal{D}}^{(t)} \rightarrow \bs{\mathcal{A}}^{(t+1)}\backslash \Phi^{(t)}_g(\bs{\mathcal{D}}^{(t)}),
\end{equation}
such that $\bs{a}^{(t)}_i\in \bs{\mathcal{A}}^{(t)}\backslash \bs{\mathcal{D}}^{(t)}$ is mapped into $\bs{a}^{(t+1)}_j \in \bs{\mathcal{A}}^{(t+1)}\backslash \Phi^{(t)}_g(\bs{\mathcal{D}}^{(t)})$ for $i, j=m+1, \dots, N$, if 
\begin{equation*}
\begin{split}
\|\bs{a}^{(t)}_i-\bs{a}^{(t+1)}_j\| \approx \|\bs{\mu}^{(t)}_1\| \ \ \text{and} \ \
(\bs{a}^{(t)}_i-\bs{a}^{(t+1)}_j)\cdot \bs{\mu}^{(t)}_1 \approx 0.
\end{split}
\end{equation*}
Figure \ref{fig:maping} also illustrates the mapping $\Phi^{(t)}_h$ and the domain $\bs{\mathcal{A}}^{(t+1)}\backslash \Phi^{(t)}_g(\bs{\mathcal{D}}^{(t)})$ for six agents. The velocity components $\bs{v}^{(t)}_i$ for $i=m+1, m+2, \dots, N$, of agents in $\bs{\mathcal{A}}^{(t)} \backslash \bs{\mathcal{D}}^{(t)}$ are defined as the displacements of position \begin{equation}\label{eqn:velocity_g}
\bs{v}^{(t)}_i=\Phi^{(t)}_h(\bs{a}^{(t)}_i)-\bs{a}^{(t)}_i.
\end{equation}

The set $\bs{\mathcal{V}}^{(t)}=\left\{\bs{v}^{(t)}_1, \bs{v}^{(t)}_2, \dots, \bs{v}^{(t)}_N\right\}$ represents all the individual velocity components at the $t$-th time step. Thus, the mapping $\Phi^{(t)}$ in (1) is defined as a union of two individual maps $\Phi^{(t)}_g$ and $\Phi^{(t)}_h$ such that
\begin{equation}\label{eqn:mapping_complete}
\Phi^{(t)}=\Phi^{(t)}_g \cup \Phi^{(t)}_h.
\end{equation}
The mean velocity, $\bs{\mu}_{\bs{\mathcal{V}}}^{(t)}=\text{mean}\left(\bs{\mathcal{V}}^{(t)}\right)$, is considered to be the velocity of the whole group at time $t$. The normalized mean speed,
\begin{equation}
\frac{\|\bs{\mu}_{\bs{\mathcal{V}}}^{(t)}\|}{\max\{\|\bs{\mu}_{\bs{\mathcal{V}}}^{(t)}\| \ ; \forall t\}},
\end{equation}
is used to detect speed changes in the group..

Coordination of agents at the $t$-th time step is measured by the polarization \cite{butail2013analysis}, 
\begin{equation}\label{eqn:polarization}
\mathcal{P}^{(t)}=\frac{1}{N} \left\| \sum_{i=1}^N\begin{bmatrix}\cos(\theta^{(t)}_i)\\\sin(\theta^{(t)}_i)\end{bmatrix} \right\|,
\end{equation}
where $\theta^{(t)}_i=\tan^{-1}\left(v_{i,2}^{(t)}/v_{i,1}^{(t)}\right)$is the orientation of the $i$-th agent moving with velocity $\bs{v}_i^{(t)}=\left[v_{i,1}^{(t)}, v_{i,2}^{(t)}\right]^{\mathrm{T}}$ computed through the aforesaid mapping of agents. Polarization ranges between 0 and 1, such that 1 identifies perfect coordination, while 0  implies no coordination.

Change in structure is measured as the number of connected components, $C^{(t)}$, of the interaction network of the agents at each time step, by representing the interaction network of agent positions as a time-varying graph. This network is made at each time step by linking the adjacent agents determined by the range search algorithm  in \cite{friedman1977algorithm}. For each agent, this algorithm searches all other agents within a given Euclidean distance 
\begin{equation}\label{eqn:epsilon}
\epsilon=\frac{2}{TN(N-1)}\sum_{t=1}^{T} \ \  \sum_{i, j=1,\dots,N} \|\bs{a}^{(t)}_i-\bs{a}^{(t)}_j\|,
\end{equation}
which is the average distance of each pair of nearest neighbors of all agents over all time steps. The normalized number of connected components \cite{tarjan1972depth}, $C^{(t)}/N\in [0, 1]$,  is used to detect changes of the structure. 

We define the set $\mathcal{X}=\{\mathcal{X}^{(t)}\}^T_{t=1}$, which contains the values generated from a convex combination of the speed, polarization, and number of clusters, all scaled between 0 and 1. That is, for $\xi_1$ and $\xi_2 \in [0, 1]$, we let
\begin{equation}
\begin{split}\label{eqn:metric_set}
\mathcal{X}^{(t)}=\xi_1 \frac{\|\bs{\mu}_{\bs{\mathcal{V}}}^{(t)}\|}{\max\{\|\bs{\mu}_{\bs{\mathcal{V}}}^{(t)}\| ; \forall t\}}+\xi_2 \mathcal{P}^{(t)}+(1-\xi_1 -\xi_2)\frac{C^{(t)}}{N},
\end{split}
\end{equation}
so that the metric distance $\Delta:\{1, \dots, T\}\times\{1, \dots, T\} \rightarrow \mathbb{R}$ is defined as
\begin{equation}\label{eqn:metric}
\Delta(t_1,t_2)=\left|\mathcal{X}^{(t_1)}-\mathcal{X}^{(t_2)}\right|.
\end{equation}
Since $\mathcal{X}^{(t)}$ is a convex combination, the relative importance of the speed, coordination and structure can be changed by varying the scalars $\xi_1$ and $\xi_2$ so that distinct group actions driven by changes in speed, coordination, and structure should be represented through distinct manifolds. In other words, two configurations, $\bs{\mathcal{A}}^{(t_1)}$ and $\bs{\mathcal{A}}^{(t_2)}$, should give a large metric distance, $\Delta(t_1,t_2)$, if they are sampled from different manifolds representing different group actions.

%%%%       Section 03:		Examples	%%%%%%%%
\section{Examples}\label{sec:examples}
Here, we evaluate the metric on three simulations of self-propelled particles \cite{vicsek1995novel}. The self-propelled particle model  \cite{vicsek1995novel} updates the position and orientation of each particle under the influence of its nearest-neighbors within a given distance. Briefly, for a nearest-neighbor set $N_i^{(t)}$, comprising all agents within a unit distance from agent $i$ (including agent $i$ itself), the Vicsek model updates the orientation $ \theta_i^{(t)} $ of the $i$-th agent at $t$-th time step as
\begin{equation}\label{eqn:vicsek1}
\theta_i^{(t+1)}=\arg{(\bs{V}_i^{(t)})}+\epsilon_i^{(t)},
\end{equation}
where $\epsilon_i^{(t)}$ is the orientation noise parameter for the $i$-th agent at the $t$-th time-step  and $\bs{V}^{(t)}_i$ is the average direction of motion of all nearest-neighbors including the agent. We assume that the noise is sampled from a uniform distribution between $\alpha_1$ and $\alpha_2$,  $\mathbb{U}[\alpha_1, \alpha_2]$. 

We augment the model by including a two-dimensional rotation matrix $R_i^t$ for the $i$-th agent at the $t$-th time-step. Such a rotation matrix can be used to change the orientation of select agents so that the  group can split and rejoin. The position $\bs{a}_i^{(t)}$ of the $i$-th agent at the $t$-th time step is therefore updated as
\begin{equation}\label{eqn:vicsek2}
\begin{split}
\bs{a}_i^{(t+1)}=\bs{a}_i^{(t)}+s_i^{(t)} R_i^{(t)} \begin{bmatrix}\cos(\theta_i^{(t)}) \\ \sin(\theta_i^{(t)})\end{bmatrix} \delta t, \\
\bs{V}_i^{(t)}=\frac{1}{|N^{(t)}_{i}|} \sum_{j\in N^{(t)}_{i}} R_j^{(t)} \begin{bmatrix}\cos(\theta_j^{(t)}) \\ \sin(\theta_j^{(t)})\end{bmatrix},
\end{split}
\end{equation} 
where $\delta t$ is a duration of a time unit and $s_i^{(t)}$ is the speed of the $i$-th agent at the $t$-th time step. The agent speed (assumed constant in the original Vicsek model) is modeled as time varying here. If $R_j^{(t)}$ is a $2\times 2$ identity matrix for all $j$ and $s_i^{(t)}$ is constant, the model (\ref{eqn:vicsek2}) specializes to the classical Vicsek model. 

In order to produce collective motion, three simulations of the augmented Vicsek model are carried out by changing $s^{(t)}_i$, $\epsilon^{(t)}_i$, and $R^{(t)}_i$ for $T$ time steps for $N$ agents in a rectangular domain of size $2L\times 2H$ with periodic boundary conditions \cite{vicsek1995novel}. Initial alignments of agents  are set to zero, and positions are sampled from a uniform distribution in a circle of radius 2 centered at $(-L+2, 0)$. The size of the time step $\delta t$ is taken to be 0.05.

The metric in equation (\ref{eqn:metric}) is computed from the positions of agents and is used to identify and characterize the distinct manifolds representing different group actions. For all the examples shown here, we set $\xi_1=1/3$ and $\xi_2=1/3$, so that the speed, coordination, and structure are weighted equally. The pairwise metric distances are then used to identify separate group actions which are then further analyzed using Isomap.

%%%%         Example 01:      Simulation of Vicsek model with switching  speed       %%%%
\subsection{A simulation of Vicsek model with switching  speed} \label{sec:example01}
To study changes in collective motion governed by differences in group speed, we simulate 50 agents ($N$) through 150 time steps ($T$) in a domain with $L=8$ and $H=5$ and impose the following speed changes,
\begin{equation}
s_i^{(t)} =
\begin{cases}
0.05+\epsilon_{s}, & \text{if }t<50 \\
0.1+\epsilon_{s}, & \text{if }50 \le t < 100\\
0.05+\epsilon_{s}, & \text{if }t\ge 100
\end{cases},
\end{equation}
where $\epsilon_{s}\in \mathbb{U}[-0.01, 0.01]$ represents individual variability in speed with respect to the group. We set the orientation noise $\epsilon_i^{(t)} \in \mathbb{U}[-0.01, 0.01]$ for all agents and choose the rotation matrix as the identity matrix for all agents. 

\begin{figure*}[hpt]
        	\centering
	\includegraphics[width=1\textwidth]{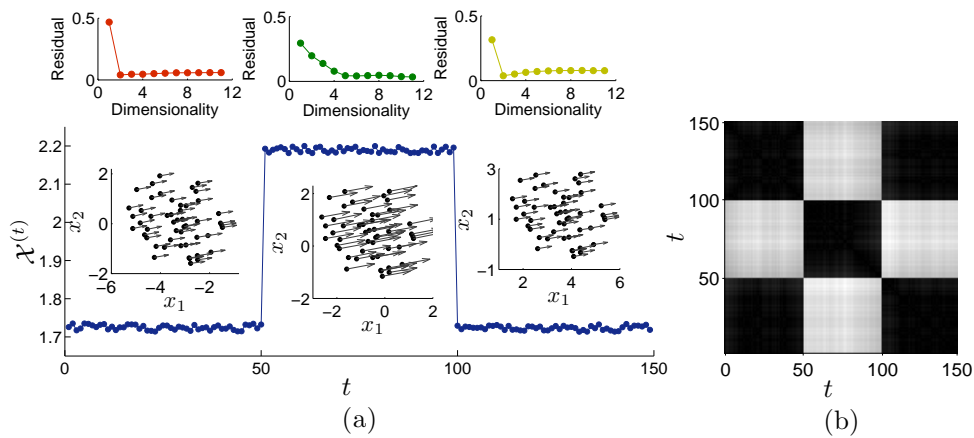}
	\vspace{-10pt}
        	\caption{(a) Distribution of $\mathcal{X}^{(t)}$ versus configurations. Snapshots  in each section show the structures and velocities of the agents. Therein, the length of the lines measures the speeds and the directions quantify the orientations of motion. Normalized Isomap residual variance versus dimensionality plots for each section are presented in the top row. For any two given configurations $\bs{\mathcal{A}}^{(t_1)}$ and $\bs{\mathcal{A}}^{(t_2)}$ where $t_1, t_2 \in 1, \dots, T$, and any arbitrary $t_0\in 1, \dots, T$, if the distance metric image (b) shows the same color at two points $(t_0, t_1)$ and $(t_0, t_2)$, then those configurations lie on the same manifold.}
	\vspace{-15pt}
\label{fig:example01}
\end{figure*}

Figure \ref{fig:example01}(a) shows that configurations for $t$ between 51 and 100 are characterized by high values of $\mathcal{X}^{(t)}$  compared to configurations through the time steps $51-100$ and $101-150$. This is evidenced from the right and left snapshots which display low speeds, while the middle snapshot that shows high speed. Residual variance plots (Fig.~\ref{fig:example01}(a) top row), obtained after running Isomap over the data, reveal that the first (red) and last (yellow) sections lie on an underlying manifold of dimensionality two. In contrast, Isomap indicates that the middle section (green) manifold has dimensionality five. Metric distances between all pairs of configurations are computed and rescaled to range between 0 and 1. Then, rescaled distances are converted to a gray color image, presented as Fig. \ref{fig:example01}(b), such that black refers to 0 and white to 1. According to the image, while configurations between $1-50$ and $101-150$ lie on the same manifold, those between $51-100$ are embedded on a different manifold.

%%%%         Example 02:      Simulation of Vicsek model with switching  direction       %%%%
\subsection{A simulation of Vicsek model with switching  coordination} \label{sec:example02}
To study changes elicited by differences in group coordination, we simulate a group of 50 agents through 150 time steps with $L$ and $H=6$ and consider the following noise on the heading of each agent in the three sections
\begin{equation}
\epsilon^{(t)}_i \in
\begin{cases}
0.01\epsilon_{c}, & \text{if }t<50 \\
0.2\epsilon_{c}, & \text{if }50 \le t < 100\\
0.01\epsilon_{c}, & \text{if }t\ge 100
\end{cases},
\end{equation}
where $\epsilon_{c}\in \mathbb{U}[-1, 1]$. The speed of agents is set to $0.05+\epsilon_{s}$ where $\epsilon_{s}\in \mathbb{U}[-0.01, 0.01]$ and the rotation matrix is chosen to be the identity matrix. 

\begin{figure*}[hpt]
        	 \centering
	\includegraphics[width=1\textwidth]{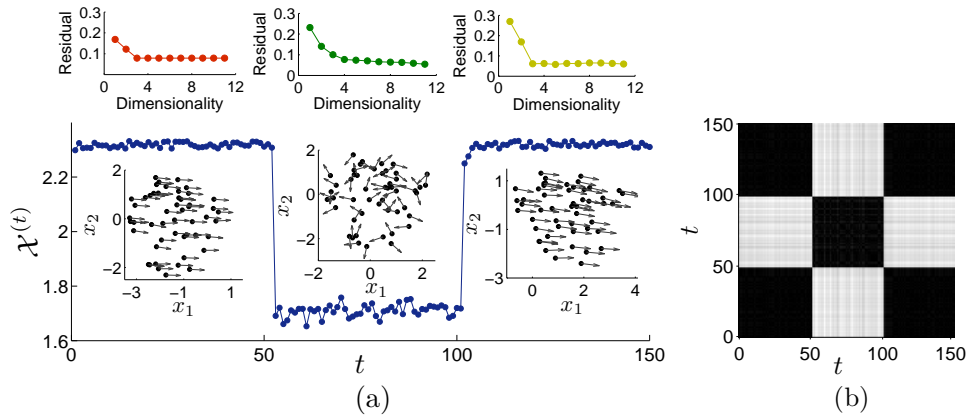}
	\caption{(a) Distribution of $\mathcal{X}^{(t)}$ versus configurations with snapshots showing the structures along with velocities of the agents. There-in, the lengths measures their speeds and arrows quantify their orientations. Normalized Isomap residual variance plots for each section (red, green, yellow) are shown in the top row. For any two given configurations $\bs{\mathcal{A}}^{(t_1)}$ and $\bs{\mathcal{A}}^{(t_2)}$ where $t_1, t_2 \in 1, \dots, T$, and any arbitrary $t_0\in 1, \dots, T$, if the distance metric image (b) gives the same color for pixels at two points $(t_0, t_1)$ and $(t_0, t_2)$, then those configurations are in the same manifold, or in different manifolds otherwise.}
	
\label{fig:example02}
\end{figure*}

Figure \ref{fig:example02}(a) demonstrates that  configurations through time steps $51-100$ are characterized by higher values of  $\mathcal{X}^{(t)}_i$ than those between time steps $1-50$ and $101-150$. Such a variation should be ascribed to the changes in the coordination of the group, as further evidenced from the snapshots presented therein. Specifically, while the right and left snapshots look similar, with agents moving with high coordination, the middle snapshot displays low coordination. The Isomap residual plot in Fig. \ref{fig:example02}(b) confirms these similarities in terms of the dimensionality, in that the first (red) and last (yellow) sections lie on an underlying manifold of dimensionality three, while the second (green) manifold has dimensionality four. As in the previous example, the metric distances between all pairs of configurations are computed and presented as an image in Fig. \ref{fig:example02}(b). Results presented therein confirm that configurations between time steps $1-50$ and $101-150$ lie on one manifold, while configurations between $51-100$ embed on a different manifold. 

%%%%         Example 03:      Simulation of Vicsek model with switching connected components       %%%%
\subsection{A simulation of Vicsek model with switching number of clusters} \label{sec:example03}
\begin{figure}
\centering
\includegraphics[width=.6\textwidth]{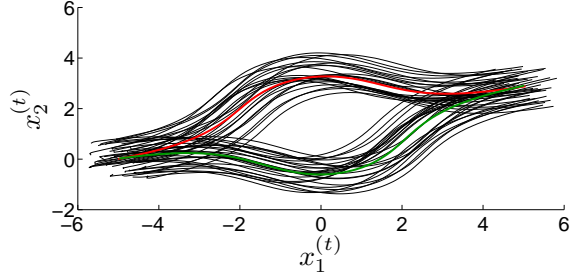}	        	
\caption{Trajectories of agents of two connected components, with red and blue curves depicting the motion of their centroids.}
\label{fig:example03_trajectory}
\end{figure}
In this example, we simulate a group of agents that initially move together, then break into two subgroups, and eventually rejoin. We use a rotation matrix to change the orientation of agents such that the two subgroups move in different directions (Fig. \ref{fig:example03_trajectory}).

We compute two-dimensional coordinates $[\bs{x}^{(1)}; \dots; \bs{x}^{(T)}]$ where $\bs{x}^{(t)}=(x_1^{(t)},x_2^{(t)})\in\mathbb{R}^2$ to obtain the movement of the centroid of the upper half of the group as shown by red in the Figure \ref{fig:example03_trajectory}. We discretize the horizontal axis into $T$ segments ranges from $-6$ to 6 by
\begin{equation}
x^{(t)}_1=6(2t-T)/T \ ; \ t=1, \dots, T.
\end{equation}
Subgroup trajectories are obtained by using difference of two sigmoid function as
\begin{equation}
x^{(t)}_2=5\left(\frac{1}{1+e^{-(\frac{T}{12}t-4)}}-\frac{1}{1+e^{-(\frac{T}{12}t-8))}}\right); t=1,\dots, T.
\end{equation}
The rotation angle between $t$-th and $(t+1)$-th time steps for half the agents is
\begin{equation}
\gamma^{(t)}=\tan\left((x_2^{(t)}-x_{2}^{(t-1)})/(x_1^{(t)}-x_{1}^{(t-1)})\right) \ ; \ t=2,\dots,T.
\end{equation} 
The rotation angles for the other half is $-\gamma^{(t)}$ for $t=2,\dots,T$. We arbitrarily split the whole group into two groups as
\begin{equation}
R^{(t)}_{i}=  
\begin{cases}
\begin{pmatrix} \cos\gamma^{(t)} & -\sin\gamma^{(t)} \\ \sin\gamma^{(t)} & \cos\gamma^{(t)} \end{pmatrix}, & \text{if } \ \ 1 \le i \le N/2, \\
\begin{pmatrix} \cos(-\gamma^{(t)}) & -\sin(-\gamma^{(t)}) \\ \sin(-\gamma^{(t)}) & \cos(-\gamma^{(t)}) \end{pmatrix}, & \text{if } \ \  N/2+1 \le i \le N.
\end{cases}
\end{equation}

This special rotation matrix is used to simulate changes in group structure for 50 agents with $L= H=6$. We select $\epsilon_i^{(t)} \in \mathbb{U}[-0.01, 0.01]$ and $s_i^{(t)}=0.05+\epsilon_{s}$ with $\epsilon_{s} \in \mathbb{U}[-0.01, 0.01]$ . 

\begin{figure*}[hpt]
        	\centering
        	\includegraphics[width=1\textwidth]{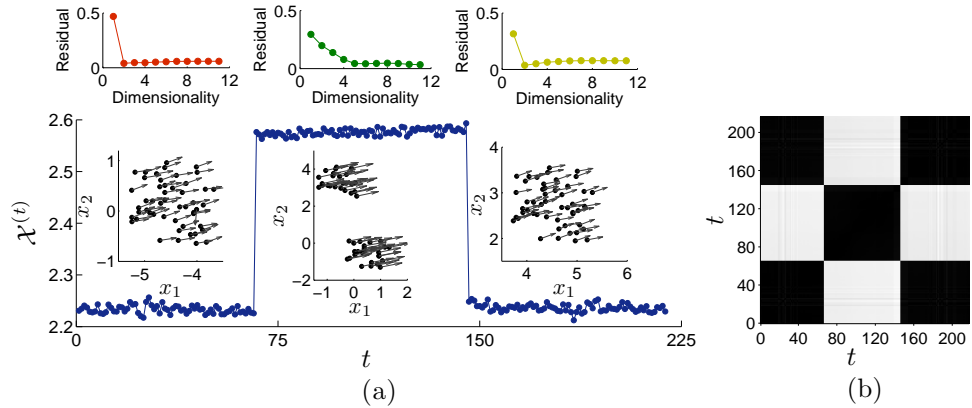}	
        	\caption{(a) Distribution of $\mathcal{X}^{(t)}$ versus configurations with snapshots giving the structures and velocities. Therein, the length of lines measures their speeds and arrows quantify the their orientations. The top row of plots shows normalized Isomap residual variance versus dimensionality for each section (red, green, yellow).  For any two given configurations $\bs{\mathcal{A}}^{(t_1)}$ and $\bs{\mathcal{A}}^{(t_2)}$ where $t_1, t_2 \in 1, \dots, T$, and any arbitrary $t_0\in 1, \dots, T$, if the distance metric image (b) gives the same color for pixels at two points $(t_0, t_1)$ and $(t_0, t_2)$, then those configurations are in the same manifold, or in different manifolds otherwise.}
	\vspace{-10pt}
\label{fig:example03}
\end{figure*}

Figure \ref{fig:example03}(a) reveals that configurations between $66-145$ result in higher $\mathcal{X}^{(t)}$ values than those in during time steps $1-65$ and $146-220$. The right and left snapshots illustrate that the group moves as one, while the middle snapshot indicates the presence of two subgroups. Isomap (Fig. \ref{fig:example03}(a)) implementation demonstrates that the first (red) and last (yellow) sections lie on an underlying manifold of dimensionality two, while the middle section (green) manifold has dimensionality four. Metric distances presented as an image (Fig. \ref{fig:example03}b) show that the configurations between $1-65$ and $146-220$ lie on one manifold and those between $66-145$ on a different manifold.

%%%%         Section 04: Discussion      %%%%
\section{Discussion} \label{sec:conclusion}
Three diverse examples of changing speed, coordination, and structure are simulated using an augmented Vicsek model to validate our novel method of identifying manifolds. Once distinct manifolds are revealed on the basis of metric distance, Isomap dimensionality reduction routine is used to compare their dimensionalities. The three examples presented here offer some validation of our method and provide useful insight into the possibility of handling alternating group actions of multi-agent systems in a fully data-driven approach. Though the Vicsek model simulates agent positions in time and records them in the same order between all time steps, our method is also robust to the unavailability of such ordering information. This is particularly important in dealing with real experimental data, where tracking individuals and preserving their identities \cite{khan2005mcmc, wu2009tracking-reconstruction, wu2009tracking-a, ristic2011metric} may be challenging.

In this study, we proposed a novel metric to characterize alternating manifolds of a group of agents. The method takes two-dimensional positions at each time step as input. Individual positions in successive configurations are mapped for each agent by a simple nearest-neighbor search in position and then velocity. While more accurate methods exist for individual tracking \cite{yilmaz2006object}, we found that this approach is $O(\log N)$ fast and robust to errors in position assignments to individual agents in successive configurations, much in the spirit of a coarse observable.

The relative importance of average speed, polarization and normalized connected components  associated with the convex representation in the proposed metric can be changed by tuning the scalers $\xi_1$ and $\xi_2$ so that the metric can effectively encounter the similarities of the group. Compared to \cite{kolpas2007coarse}, where average nearest-neighbor distance is used to analyze group motion, the combined metric proposed here captures change across multiple group motions. Different than \cite{kolpas2007coarse}, the connected-component weighting is able to also identify differences in group structure in terms of group splits.

The proposed metric is able to isolate distinct phases of collective motion in three separate scenarios. The separation of these phases as sub-manifolds in a low-dimensional space was further confirmed by the Isomap algorithm, which measures a different dimensionality for each of these phases. Further based on metric distances it was possible to group similar actions together. This has potential application in using the proposed metric to train specific events of interest and process long videos to highlight the same. 

Multi-agent groups often change their behavior due to natural perturbations \cite{carere2009aerial, couzin2002collective, helbing2000simulating, neill1974experiments, vicsek2012collective}. We posit that the entire evolution can be described by a collection of sub-manifolds, each representing similar group actions. Identifying such sub-manifolds can offer important insight on the analysis the collective motion and form an important step before more detailed investigations are conducted. While the artificial separation of phases was distinct in the simulated examples presented here, future work will evaluate the performance and sensitivity of this metric to gradual changes over a sequence of time-steps. Future work will also focus on application of this metric to analyze videos from experiments on collective behavior.

\section*{Acknowledgments}
Kelum Gajamannage and Erik M. Bollt have been supported by the National Science Foundation under grant no. CMMI-1129859.
Sachit Butail and Maurizio Porfiri have been supported by the National Science Foundation under grant no. CMMI-1129820.


\begin{thebibliography}{10}

\bibitem{abaid2012topological}
N.~Abaid, E. Bollt, M.~Porfiri,
\newblock Topological analysis of complexity in multiagent systems.
\newblock Physical Review E, 85(4):041907, (2012).

\bibitem{butail2013analysis}
S.~Butail, E.M. Bollt, M.~Porfiri,
\newblock Analysis and classification of collective behavior using generative
  modeling and nonlinear manifold learning.
\newblock Journal of Theoretical Biology, 336:185--199, (2013).

\bibitem{carere2009aerial}
C.~Carere, S.~Montanino, F.~Moreschini, F.~Zoratto, F.~Chiarotti, D.~Santucci, E.~Alleva,
\newblock Aerial flocking patterns of wintering starlings, Sturnus
  vulgaris, under different predation risk.
\newblock Animal Behaviour, 77(1):101--107, (2009).

\bibitem{couzin2009collective}
I.D. Couzin,
\newblock Collective cognition in animal groups.
\newblock Trends in Cognitive Sciences, 13(1):36--43, (2009).

\bibitem{couzin2002collective}
I.D. Couzin, J.~Krause, R.~James, G.D. Ruxton, N.R. Franks,
\newblock Collective memory and spatial sorting in animal groups.
\newblock Journal of Theoretical Biology, 218(1):1--11, (2002).

\bibitem{delellis2014collective}
P.~DeLellis, G.~Polverino, G.~Ustuner, N.~Abaid, S.~Macr{\`\i}, E.M. Bollt, M.~Porfiri,
\newblock Collective behaviour across animal species.
\newblock Scientific Reports, 4, (2014).

\bibitem{delellis2013topological}
P.~DeLellis, M.~Porfiri, E.M. Bollt,
\newblock Topological analysis of group fragmentation in multiagent systems.
\newblock Physical Review E, 87(2):022818, (2013).

\bibitem{deneubourg1989collective}
J.~Deneubourg, S.~Goss,
\newblock Collective patterns and decision-making.
\newblock Ethology Ecology \& Evolution, 1(4):295--311, (1989).

\bibitem{friedman1977algorithm}
J.~H. Friedman, J.~L. Bentley, R.~A. Finkel,
\newblock An algorithm for finding best matches in logarithmic expected time.
\newblock ACM Transactions on Mathematical Software (TOMS), 3(3):209--226, (1977).

\bibitem{gajamannage2015dimensionality}
K.~Gajamannage, S.~Butail, M.~Porfiri, E.M. Bollt,
\newblock Dimensionality reduction of collective motion by principal manifolds.
\newblock  Physica D: Nonlinear Phenomena, 291:62--73, (2015).

\bibitem{helbing2000simulating}
D.~Helbing, I.~Farkas, T.~Vicsek,
\newblock Simulating dynamical features of escape panic.
\newblock Nature, 407(6803):487--490, (2000).

\bibitem{helbing2002simulation}
D.~Helbing, I.J. Farkas, P.~Molnar, T.~Vicsek,
\newblock Simulation of pedestrian crowds in normal and evacuation situations.
\newblock Pedestrian and Evacuation Dynamics, 21:21--58, (2002).

\bibitem{khan2005mcmc}
Z.~Khan, T.~Balch, and F.~Dellaert,
\newblock Mcmc-based particle filtering for tracking a variable number of
  interacting targets.
\newblock Pattern Analysis and Machine Intelligence, IEEE Transactions  on, 27(11):1805--1819, (2005).

\bibitem{kolpas2007coarse}
A.~Kolpas, J.~Moehlis, I.G. Kevrekidis,
\newblock Coarse-grained analysis of stochasticity-induced switching between
  collective motion states.
\newblock Proceedings of the National Academy of Sciences, 104(14):5931--5935, (2007).

\bibitem{lima1990behavioral}
S.L. Lima, L.M. Dill,
\newblock Behavioral decisions made under the risk of predation: a review and
  prospectus.
\newblock Canadian Journal of Zoology, 68(4):619--640, (1990).

\bibitem{neill1974experiments}
S.R. Neill, J.M. Cullen,
\newblock Experiments on whether schooling by their prey affects the hunting
  behaviour of cephalopods and fish predators.
\newblock Journal of Zoology, 172(4):549--569, (1974).

\bibitem{parrish1991predators}
J.K. Parrish,
\newblock Do predators shape fish schools: interactions between predators and
  their schooling prey.
\newblock Netherlands Journal of Zoology, 42(2):358--370, (1991).

\bibitem{ristic2011metric}
B.~Ristic, B.~Vo, D.~Clark, B.~Vo,
\newblock A metric for performance evaluation of multi-target tracking
  algorithms.
\newblock IEEE Transactions on Signal Processing, 59(7):3452--3457, (2011).

\bibitem{shiwakoti2011animal}
N.~Shiwakoti, M.~Sarvi, G.~Rose, M.~Burd,
\newblock Animal dynamics based approach for modeling pedestrian crowd egress
  under panic conditions.
\newblock Transportation Research part B: methodological,
  45(9):1433--1449, (2011).

\bibitem{tarjan1972depth}
R.~Tarjan,
\newblock Depth-first search and linear graph algorithms.
\newblock SIAM Journal on Computing, 1(2):146--160, (1972).

\bibitem{tenenbaum2000global}
J.B. Tenenbaum, V.~De~Silva, J.C. Langford,
\newblock A global geometric framework for nonlinear dimensionality reduction.
\newblock Science, 290(5500):2319--2323, (2000).

\bibitem{vabo1997individual}
R.~Vabo, L.~Nottestad,
\newblock An individual based model of fish school reactions: predicting
  antipredator behaviour as observed in nature.
\newblock Fisheries Oceanography, 6(3):155--171, (1997).

\bibitem{vicsek1995novel}
T.~Vicsek, A.~Czir{\'o}k, E.~Ben-Jacob, I.~Cohen, O.~Shochet,
\newblock Novel type of phase transition in a system of self-driven particles.
\newblock Physical Review Letters, 75(6):1226, (1995).

\bibitem{wu2009tracking-a}
Z.~Wu, N.I. Hristov, T.L. Hedrick, T.H. Kunz, M.~Betke,
\newblock Tracking a large number of objects from multiple views.
\newblock In {\em Proceedings of the IEEE International Conference on Computer
  Vision}, pages 1546--1553. IEEE, (2009).

\bibitem{wu2009tracking-reconstruction}
Z.~Wu, N.I. Hristov, T.H. Kunz, M.~Betke,
\newblock Tracking-reconstruction or reconstruction-tracking? comparison of two
  multiple hypothesis tracking approaches to interpret 3d object motion from
  several camera views.
\newblock In {\em Proceedings of the Workshop on Motion and Video Computing},
  pages 1--8. IEEE, (2009).

\bibitem{yilmaz2006object}
A.~Yilmaz, O.~Javed, M.~Shah,
\newblock Object tracking: A survey.
\newblock ACM Computing Surveys (CSUR), 38(4):13, (2006).

\bibitem{vicsek2012collective}
T.~Vicsek, A.~Zafeiris,
\newblock Collective motion.
\newblock Physics Reports, 517(3):71--140, (2012).

\end{thebibliography}
\end{document}